\documentclass[11pt,reqno]{amsart}
\usepackage[margin=1in]{geometry}

\usepackage{amssymb,amsmath,amsthm}
\usepackage{microtype}
\usepackage{mathtools}
\usepackage{enumitem}

\usepackage[colorlinks=true, pdfstartview=FitH, linkcolor=blue, citecolor=blue, urlcolor=blue]{hyperref}





\renewcommand{\le}{\leqslant}

\renewcommand{\ge}{\geqslant}

\newcommand{\norm}[1]{\left\lVert#1\right\rVert_2}
\newcommand{\inner}[2]{\left\langle#1, #2\right\rangle}

\newcommand{\N}{\mathbb{N}}
\newcommand{\R}{\mathbb{R}} 
\newcommand{\Z}{\mathbb{Z}}
\newcommand{\F}{\mathbb{F}}
\newcommand{\ex}{\operatorname{ex}}
\newcommand{\vol}{\operatorname{vol}}

\newtheorem{theorem}{Theorem}[section]

\newtheorem{lemma}[theorem]{Lemma}

\theoremstyle{definition}

\numberwithin{theorem}{section}

\allowdisplaybreaks

\begin{document}

\title{On the number of triangles in $K_{1,2,2}$-free graphs}

\author{Ritesh Goenka}

\address{Mathematical Institute \\ University of Oxford \\ Oxford OX2 6GG, United Kingdom }
\email{goenka@maths.ox.ac.uk}


\begin{abstract}
    For every sufficiently large $n \in \mathbb{N}$, we construct an $n$-vertex tripartite graph with at least $n^{2.0002}$ triangles in which no vertex neighbourhood contains a $4$-cycle, thereby giving an improved lower bound on the generalised Tur\'an number $\operatorname{ex}(n, K_3, K_{1,2,2})$. This answers a question of Kostochka, Mubayi, and Verstra\"ete on the Tur\'an numbers of expansions of even wheels and also disproves Tsai's induced double Tur\'an conjecture. Our construction labels vertices by vectors in the union of a lattice and one of its translates so that the three vertex labels on every triangle sum to an element of a Sidon set, thereby forbidding a $4$-wheel. The proof is due to GPT~6.2 Astra.
\end{abstract}

\maketitle

\section{Introduction}

Classical Tur\'an problems ask for the maximum number of edges in an $n$-vertex graph if it contains no copy of a given graph $H$. The idea of maximising the count of other subgraphs goes back at least to Zykov's work~\cite{Zyk} on cliques in $K_r$-free graphs (see the survey of Gerbner and Palmer~\cite{GP}). Alon and Shikhelman~\cite{AS} initiated the systematic study of the generalised Tur\'an number $\ex(n,T,H)$, which denotes the maximum number of copies of $T$ in an $H$-free graph on $n$ vertices. For $T = K_3$, they determined the asymptotics of $\ex(n,K_3,H)$ whenever $\chi(H) \ge 4$. However, for $\chi(H) = 3$, even the correct power of $n$ can be difficult to determine.

Perhaps the simplest case in which the known upper and lower bounds differ by a polynomial factor is $\ex(n,K_3,K_{1,2,2})$. The graph $W_4 = K_{1,2,2}$ is the $4$-wheel, obtained by joining a vertex to all four vertices of a cycle of length four. Adding a maximal matching inside one part of a balanced complete bipartite graph gives a $W_4$-free graph with $\Omega(n^2)$ triangles. In the other direction, Mubayi and Mukherjee~\cite{MM} proved an upper bound of $o(n^{5/2})$, which was also proved independently by Methuku, Gr\'osz, and Tompkins in an unpublished work (see~\cite{MM}) and by Balogh (see his lecture note~\cite{Bal}). Calbet and the author~\cite{CG} subsequently improved this upper bound to $O(n^{5/2}/\sqrt{\log n})$. The problem of narrowing the gap between the bounds has been raised by several authors~\cite{BJL,MM,MV}. We give a polynomial improvement to the lower bound.

\begin{theorem}
\label{thm:main}
    For every sufficiently large $n \in \N$, there exists a tripartite $W_4$-free graph on $n$ vertices with at least $n^{2.0002}$ triangles. Thus $\ex(n, K_3, W_4) = \Omega(n^{2.0002})$.
\end{theorem}

Despite the improvement in the lower bound, there remains a polynomial gap in the bounds for $\ex(n,K_3,W_4)$. Theorem~\ref{thm:main} answers a question of Kostochka, Mubayi and Verstra\"ete~\cite[Question~2]{KMV}, who asked whether the expansion\footnote{Let $G$ be a graph with vertex set $V(G)$ and edge set $E(G)$. Its expansion is the $3$-uniform hypergraph $G^+$ with edge set $\{e \cup \{v_e\} : e \in E(G)\}$, where $(v_e)_{e \in E(G)}$ are distinct vertices not in $V(G)$.} of every even wheel has Tur\'an number $O(n^2)$. Indeed, the triangles in any graph given by Theorem~\ref{thm:main} form a $3$-uniform hypergraph with no copy of the expansion of $W_4$, since such a copy would give a copy of $W_4$ in the original graph.

Theorem~\ref{thm:main} also disproves Tsai's conjecture on induced double Tur\'an numbers~\cite[Conjecture~A]{Tsai}. For any fixed forest $F$, the conjecture predicts that $n$ induced subgraphs of an $n$-vertex graph have $O_F(n^2)$ edges in total whenever their pairwise intersections are $F$-free. Now let $G$ be an $n$-vertex graph given by Theorem~\ref{thm:main} and consider the $n$ induced subgraphs $G[N(v)]$ for $v \in V(G)$. Suppose that the common neighbourhood of two distinct vertices $u$ and $v$ contains a path $xyz$ on three vertices. Then $G$ contains a $4$-wheel with central vertex $y$ and rim vertices $u, x, v, z$ in cyclic order, a contradiction. Thus the pairwise intersections are $F$-free for every forest $F$ containing a path on three vertices. Yet each triangle of $G$ contributes one edge to the neighbourhood graph of each of its three vertices, so the total number of edges in the neighbourhood graphs is superquadratic.

The first step of our proof involves constructing a lattice $\Lambda$ and a translate $\Gamma$ of this lattice such that the set $\Gamma_{\min}$ of shortest vectors in $\Gamma$ is \emph{Sidon}. We then consider a complete tripartite graph with vertices in two parts labelled by points of $\Lambda$ and vertices in the third part labelled by points of $\Gamma$, so that the sum of the vertex labels on every triangle lies in $\Gamma$. We further assign real labels to the edges so that for every triangle, the sum of its edge labels equals the squared norm of the sum of its vertex labels. Among the triangles whose vertex labels sum to a vector in $\Gamma_{\min}$, choose a most frequent triple of edge labels, read in a fixed order of the three pairs of parts. Let $G$ be the graph obtained by retaining the edges with the corresponding label between each pair of parts. The identity above ensures that the vertex labels of every triangle in $G$ sum to a vector in $\Gamma_{\min}$. A $4$-wheel $W$ in $G$ would give two distinct pairs of triangle sums in $\Gamma_{\min}$ with the same total, since each pair of opposite triangles in $W$ counts the central vertex twice and every rim vertex once. This contradicts the Sidon property of $\Gamma_{\min}$, so $G$ is $W_4$-free.

The remainder of this paper is organised as follows. In Section~\ref{sec:sidon}, we give the underlying Sidon set construction. In Section~\ref{sec:proof}, we use it to prove our main result.

\vspace{4pt}
\paragraph{\textbf{Notation.}} For $A, B \ge 0$, we write $A \lesssim B$, or equivalently $B \gtrsim A$, if $A \le CB$ for some absolute constant $C > 0$ throughout the stated parameter range. We write $A \asymp B$ if both $A \lesssim B$ and $B \lesssim A$.

\section*{Acknowledgements}

The author thanks Jacques Verstra\"ete for presenting Ray Tsai's induced double Tur\'an conjecture at the \emph{Innovations and Challenges in Extremal Combinatorics} workshop, held at the Institute for Mathematical Sciences, National University of Singapore. The author is supported by a joint Clarendon Fund and Exeter College SKP scholarship.

\section{A Sidon set of shortest vectors in a lattice translate}
\label{sec:sidon}

We begin with some definitions. Let $d \in \N$. A set $S \subseteq \R^d$ is said to be \emph{Sidon} if, whenever $x, y, z, w \in S$ satisfy $x + y = z + w$, the multisets $\{x,y\}$ and $\{z,w\}$ are equal.

For a full-rank lattice $\Lambda \subset \R^d$, we denote its covolume by $\det \Lambda$. Its \emph{covering radius} is the least $r > 0$ such that every point of $\R^d$ lies within distance $r$ of some point in $\Lambda$. For a translate $\Gamma = T + \Lambda$ of $\Lambda$, where $T \in \R^d$, we denote by $\Gamma_{\min}$ the set of shortest (in the Euclidean norm) vectors in $\Gamma$.

Let $\F_q$ be a finite field and let $k \in \N$. A \emph{linear code} of length $k$ over $\F_q$ is a subspace $\mathcal{C} \subseteq \F_q^k$, whose elements are called \emph{codewords}. We call $\mathcal{C}$ \emph{intersecting} if any two nonzero codewords in $\mathcal{C}$ have intersecting supports, that is, if they are both nonzero in some common coordinate.

We will require the following lemma, which corresponds to the case $q = 3$ and $t = 0$ of Sagalovich's random-coding bound~\cite[\S4, Eq.~(13)]{Sagalovich}. We include a short proof for completeness.

\begin{lemma}
\label{lem:code}
    Let $k \ge 4$ be an integer and let $r = \lfloor k \log_3 (3/\sqrt5) \rfloor$. Then there exists an $r$-dimensional intersecting linear code $\mathcal{C}$ in $\F_3^k$. 
\end{lemma}

\begin{proof}
    Choose a linear map $\varphi \colon \F_3^r \to \F_3^k$ uniformly at random. For fixed linearly independent vectors $u, v \in \F_3^r$, the coordinate pairs $(\varphi(u)_i,\varphi(v)_i)$ for $i \in [k]$ are independent and uniform on $\F_3^2$. Exactly five of the nine pairs in $\F_3^2$ have at least one zero entry, so $\varphi(u)$ and $\varphi(v)$ have disjoint supports with probability $(5/9)^k$. Moreover, each fixed nonzero vector lies in $\ker\varphi$ with probability $3^{-k}$. By the union bound, the probability that $\varphi$ has a nontrivial kernel or maps some linearly independent pair of vectors to vectors with disjoint supports is at most
    \begin{equation*}
        3^{r-k} + \frac{1}{2} 3^{2r} \left(\frac{5}{9}\right)^k \le \frac{1}{5^{k/2}} + \frac{1}{2} < 1.
    \end{equation*}
    Thus, there exists a linear map $\varphi_0 : \F_3^r \rightarrow \F_3^k$ with trivial kernel such that for any two linearly independent vectors $u, v \in \F_3^r$, $\varphi_0(u)$ and $\varphi_0(v)$ have intersecting supports. Further, it follows from the linearity of $\varphi_0$ that for any two linearly dependent nonzero vectors $u, v \in \F_3^r$, $\varphi_0(u)$ and $\varphi_0(v)$ have the same support. Hence, $\mathcal{C} \coloneqq \varphi_0(\F_3^r)$ is an intersecting $r$-dimensional code in $\F_3^k$.
\end{proof}

\begin{lemma}
\label{lem:lattice}
    For every integer $k \ge 4$, there exist a full-rank lattice $\Lambda \subset \R^{2k}$ and a translate $\Gamma$ of $\Lambda$ with the following properties.
    \begin{enumerate}
        \item $\det \Lambda \asymp (\sqrt{15}/2)^k$.
        \item $\Gamma_{\min} = \{s \in \Gamma: \norm{s}^2 = k/3\}$.
        \item $\Gamma_{\min}$ is Sidon and satisfies $|\Gamma_{\min}| \asymp (3/\sqrt5)^k$.
        \item The covering radius of $\Lambda$ is at most $h \coloneqq 3\sqrt{k}$.
        \item For all $x,y \in \Lambda \cup \Gamma$, we have $18 \inner{x}{y} \in \Z$.
    \end{enumerate}
\end{lemma}

\begin{proof}
    We first construct a planar lattice translate with three shortest vectors. Let $u, v \in \R^2$ be unit vectors with $\inner{u}{v} = 1/2$, and let $L = \Z u + \Z v$. Set $t = -(u+v)/3$, so that $-t$ is the centre of the equilateral triangle with vertices $0,u,v$. It is easy to check that the shortest vectors in $t+L$ are $t$, $t+u$, and $t+v$, each of squared Euclidean norm $1/3$.

    Now define an additive group homomorphism $\pi \colon L \to \F_3$ by $\pi(u) = 1$ and $\pi(v) = 2$. Let $\mathcal{C}$ be the code in $\F_3^k$ obtained from Lemma~\ref{lem:code}, and let $\pi^k \colon L^k \to \F_3^k$ act coordinatewise. Define
    \begin{equation*}
        \Lambda \coloneqq (\pi^k)^{-1}(\mathcal{C}) \subseteq L^k, \qquad T \coloneqq (t,\ldots,t), \qquad \Gamma \coloneqq T + \Lambda.
    \end{equation*}
    Each vector in $T + L^k = (t+L)^k$ has squared Euclidean norm at least $k/3$, with equality precisely when each of its planar coordinates is one of $t$, $t+u$, $t+v$. Note that the map $w \mapsto \pi(w-t)$ sends $t$, $t+u$, $t+v$ to $0$, $1$, $2$, respectively. Therefore, the map $\varphi \colon (T + L^k)_{\min} \to \F_3^k$ defined by $\varphi(p) = \pi^k(p - T)$ is a bijection. Since $T \in \Gamma$ attains the minimum squared norm $k/3$, the definition of $\Gamma$ implies that $\varphi$ restricts to a bijection from $\Gamma_{\min}$ to $\mathcal{C}$. Hence, $|\Gamma_{\min}| = |\mathcal{C}| = 3^r \asymp (3/\sqrt{5})^k$.

    To prove that $\Gamma_{\min}$ is Sidon, suppose that $p_1+p_2 = q_1+q_2$ for some $p_1,p_2,q_1,q_2 \in \Gamma_{\min}$. In each planar coordinate, the four entries belong to $\{t,t+u,t+v\}$. These three points are affinely independent, so the sum of two of them determines the two summands up to order. Hence, in every coordinate, $p_1$ agrees with either $q_1$ or $q_2$. It follows that the two codewords
    \begin{equation*}
        \varphi(p_1) - \varphi(q_1) \qquad \text{and} \qquad \varphi(p_1) - \varphi(q_2)
    \end{equation*}
    have disjoint supports. Since $\mathcal{C}$ is intersecting, at least one of these codewords is zero, which implies $p_1 = q_1$ or $p_1 = q_2$. The equality of the sums now gives equality of the remaining summands.

    Since $\pi^k$ is surjective, $\Lambda$ has index $|\F_3^k|/|\mathcal{C}| = 3^{k-r}$ in $L^k$. As $\det L = \sqrt{3}/2$, we obtain
    \begin{equation*}
        \det\Lambda
        = 3^{k-r}\left(\frac{\sqrt{3}}{2}\right)^k
        \asymp \left(\frac{\sqrt{15}}{2}\right)^k.
    \end{equation*}
    Note that $(3L)^k \subseteq \Lambda$ since $\pi^k((3L)^k) = \{0\} \subseteq \mathcal{C}$. To bound the covering radius of $\Lambda$, write a point of $\R^2$ as $au+bv$ and round $a$ and $b$ to the nearest multiples of $3$. Each coefficient changes by at most $3/2$, so the resulting point of $3L$ lies within distance $3$ of the original point. Applying this in each of the $k$ coordinates gives a point of $(3L)^k \subseteq \Lambda$ within distance $3 \sqrt{k}$ of any point of $\R^{2k}$.

    Finally, we have $\inner{u}{u} = \inner{v}{v} = 1$ and $\inner{u}{v} = 1/2$, so $2 \inner{p}{q} \in \Z$ for all $p,q \in L^k$. Since both $\Lambda$ and $3\Gamma$ lie in $L^k$, it follows that $18 \inner{x}{y} \in \Z$ for all $x, y \in \Lambda \cup \Gamma$.
\end{proof}

\section{Proof of the main result}
\label{sec:proof}

Let $d \in \N$. For $R > 0$, let $B_R^d$ denote the closed ball of radius $R$ centred at the origin. We shall omit the superscript when the ambient dimension is clear from context. Further, we denote the Lebesgue measure on $\R^d$ by $\vol_d(\cdot)$. 

\begin{lemma}
\label{lem:main}
    Let $k \ge 4$ be an integer and $R \ge 1$. Let $\Lambda, \Gamma \subset \R^{2k}$ be a lattice and one of its translates obtained from Lemma~\ref{lem:lattice}. Let $X$, $Y$, $Z$ be disjoint vertex classes, labelled by $\Lambda\cap B_R$, $\Lambda\cap B_R$, $\Gamma \cap B_R$, respectively. Let $\Delta$ denote the number of triples $(x,y,z) \in X \times Y \times Z$ with $x+y+z \in \Gamma_{\min}$. Then there is a tripartite $W_4$-free graph on these classes with at least $\Delta/(16 R)^6$ triangles.
\end{lemma}

\begin{proof}
    We may assume $\Delta > 0$ since otherwise the result follows trivially. Let $H$ be the complete tripartite graph with parts $X$, $Y$, and $Z$. For each edge $(w_1,w_2)$ in $H$, assign it the label
    \begin{equation*}
        \ell(w_1,w_2) = \frac{1}{2} \bigl(\norm{w_1}^2 + \norm{w_2}^2\bigr) + 2\inner{w_1}{w_2}.
    \end{equation*}
    The labels have absolute value at most $3R^2$ and belong to $\tfrac{1}{36} \Z$ by Lemma~\ref{lem:lattice}~(e). Thus the total number of labels is at most $(16 R)^2$. For each triple $(x, y, z) \in X \times Y \times Z$, record the ordered triple of edge labels $(\ell(x,y), \ell(y,z), \ell(z,x))$. Now choose the most frequent ordered triple $(\ell_1, \ell_2, \ell_3)$ of edge labels among the $\Delta$ triples from the statement of the lemma, and retain only the edges with the chosen label for their class pair. This gives us a graph $G$ with at least $\Delta/(16 R)^6$ triangles.

    Let $(x_0, y_0, z_0)$ be one of the $\Delta$ triples with edge labels $(\ell_1, \ell_2, \ell_3)$. Then we have
    \begin{equation*}
        \ell_1 + \ell_2 + \ell_3 = \ell(x_0,y_0) + \ell(y_0,z_0) + \ell(z_0,x_0) = \norm{x_0+y_0+z_0}^2 = k/3,
    \end{equation*}
    where the last equality follows from Lemma~\ref{lem:lattice}~(b). Thus, for any triangle $(x, y, z)$ in $G$, we have
    \begin{equation*}
        \norm{x + y + z}^2 = \ell(x,y) + \ell(y,z) + \ell(z,x) = \ell_1 + \ell_2 + \ell_3 = k/3.
    \end{equation*}
    Moreover, $x + y + z \in \Gamma$ since $x, y \in \Lambda$ and $z \in \Gamma$. It follows from Lemma~\ref{lem:lattice}~(b) that $x + y + z \in \Gamma_{\min}$.
    
    Now suppose $G$ contains a $4$-wheel $W$. Let $a$ be the label of the central vertex of $W$ and let $b_1, b_2$ and $c_1, c_2$ be labels of the rim vertices in the remaining two classes. Then its four triangle sums $s_{ij} = a + b_i + c_j$, $i, j \in [2]$, all lie in $\Gamma_{\min}$ and satisfy $s_{11} + s_{22} = s_{12} + s_{21}$. Since $\Gamma_{\min}$ is Sidon, either $s_{11} = s_{12}$ or $s_{11} = s_{21}$, implying $c_1 = c_2$ or $b_1 = b_2$ respectively. This is a contradiction since the two rim vertices in each class have distinct labels.
\end{proof}

\begin{proof}[Proof of Theorem~\ref{thm:main}]
    Let $k \ge 4$ be an integer to be chosen later and let $d = 2k$. Let $h = 3 \sqrt{k}$ be the covering radius bound from Lemma~\ref{lem:lattice}~(d). Let $M = 10^5$ and $R = (M-3) \sqrt{k}$, so that $R + h = M \sqrt{k}$. We apply Lemma~\ref{lem:main} with parameters $k$ and $R$ as above.
    
    We first bound the sizes of the vertex classes in Lemma~\ref{lem:main}. We take Voronoi cells of the vertex labels with respect to $\Lambda$ for $X$ and $Y$, and with respect to $\Gamma$ for $Z$. Each such cell lies in $B_{R+h}$ and has volume $\det \Lambda$. For a fixed vertex class, the corresponding Voronoi cells overlap only on sets of volume zero, so the class has size at most
    \begin{equation*}
        N \coloneqq \frac{\vol_d(B_{R+h})}{\det\Lambda}.
    \end{equation*}

    We will now estimate the quantity $\Delta$ in Lemma~\ref{lem:main}. We say that a triple $(x, y, z) \in (\R^d)^3$ is \emph{good} if $x, y \in \Lambda \cap B_R$, $z \in \Gamma \cap B_R$, and $x + y + z \in \Gamma_{\min}$. We will first estimate the number of good triples with a given sum $s \in \Gamma_{\min}$. Set $R' = R - 7 \sqrt{k}$, so that $R' + 2h + \norm{s} \le R$. Let
    \begin{equation*}
        P = \{(a,b) \in B_{R'} \times B_{R'} : a+b \in B_{R'}\} \qquad \text{and} \qquad Q = \{(x, y) \in \Lambda \times \Lambda: (V_x \times V_y) \cap P \ne \emptyset\},
    \end{equation*}
    where $V_w$ denotes the Voronoi cell of a point $w$ in $\Lambda$. Suppose that $(x, y) \in Q$. Then, we obtain
    \begin{equation*}
        \norm{x} \le R' + h \le R, \qquad \norm{y} \le R' + h \le R, \qquad \text{and} \qquad \norm{x+y} \le R' + 2h,
    \end{equation*}
    since each Voronoi cell lies within distance $h$ of its centre. Define $z \coloneqq s - x - y$. Then $z \in \Gamma$ and $\norm{z} \le R' + 2h + \norm{s} \le R$. Thus $x, y \in \Lambda \cap B_R$ and $z \in \Gamma \cap B_R$ with $x + y + z = s$. So each pair $(x, y) \in Q$ can be extended to a good triple $(x, y, z)$ with sum $s$. The products $(V_x \times V_y)$, as $(x,y)$ varies over $\Lambda \times \Lambda$, tile the ambient space $\R^{2d}$ and each of them has volume $(\det\Lambda)^2$. Hence, at least $\vol_{2d}(P)/(\det\Lambda)^2$ of them meet $P$. Thus the number $\Delta$ of good triples satisfies
    \begin{equation}
    \label{eqn:delta}
        \Delta \ge \sum_{s \in \Gamma_{\min}} |Q| \ge |\Gamma_{\min}|\,\frac{\vol_{2d}(P)}{(\det\Lambda)^2}.
    \end{equation}
    We will now estimate $\vol_{2d}(P)$. For fixed $a \in B_{R'}$, the set of $b \in B_{R'}$ such that $(a,b) \in P$ is the lens $B_{R'} \cap (B_{R'}-a)$. The section of the lens halfway between its centres has radius at least $\sqrt{3}\,R'/2$, and the lens contains the symmetric double cone over this section of total height $R'$. Comparing this cone with the cylinder of radius $R'$ and height $2R'$ containing $B_{R'}$ gives
    \begin{equation*}
        \vol_d \bigl(B_{R'} \cap(B_{R'}-a)\bigr) \ge \frac{1}{2d}\left(\frac{\sqrt{3}}{2}\right)^{d-1}\vol_d(B_{R'}) \ge \frac{1}{4k}\left(\frac34\right)^k\vol_d(B_{R'}).
    \end{equation*}
    Integrating over $a \in B_{R'}$, we obtain
    \begin{equation*}
        \vol_{2d}(P) \ge \frac{1}{4k}\left(\frac{3}{4}\right)^k\vol_d(B_{R'})^2,
    \end{equation*}
    which when used in \eqref{eqn:delta} gives
    \begin{equation*}
        \Delta \ge \frac{|\Gamma_{\min}|}{4k} \left(\frac34\right)^k \left(\frac{\vol_d(B_{R'})}{\det\Lambda}\right)^2.
    \end{equation*}

    By Lemma~\ref{lem:main}, there exists a tripartite $W_4$-free graph $G$ on at most $3N$ vertices with $t(G) \gtrsim \Delta/R^6$, where $t(G)$ denotes the number of triangles in $G$. Since $R'/(R+h) = 1-10/M$, volume scaling gives $\vol_d(B_{R'})/\det\Lambda = N(1-10/M)^{2k}$. Substituting this into our lower bound on $\Delta$ and using $R \asymp \sqrt{k}$ and the estimate $|\Gamma_{\min}| \asymp (3/\sqrt{5})^k$ from Lemma~\ref{lem:lattice}~(c), we obtain
	\begin{equation*}
		t(G) \gtrsim \frac{N^2}{k^4}\,\rho^k, \qquad \text{where } \rho \coloneqq \frac{9}{4\sqrt{5}}\left(1-\frac{10}{M}\right)^4.
	\end{equation*}

	It remains to choose $k$ in terms of $n$. By the volume formula for a ball in $\R^{2k}$, Stirling's formula, and the determinant estimate from Lemma~\ref{lem:lattice}~(a), we have
	\begin{equation*}
		N = \frac{(\pi M^2 k)^k}{k!\det\Lambda} \asymp \frac{K^k}{\sqrt{k}}, \qquad \text{where } K \coloneqq \frac{2\pi e M^2}{\sqrt{15}}.
	\end{equation*}
    Suppose $n \ge K^4$. Take $k = \lfloor \log_K n \rfloor$, so that $k \ge 4$ and $n/K < K^k \le n$. It follows that $N \asymp n/\sqrt{k}$ and $\rho^k \asymp n^{\log_K \rho}$. Since $k \asymp \log n$, we have $3N \le n$ for all sufficiently large $n$, and we may add isolated vertices to $G$ to obtain a tripartite $W_4$-free graph on $n$ vertices with
	\begin{equation*}
		t(G) \gtrsim \frac{n^{2+\log_K \rho}}{(\log n)^5}.
	\end{equation*}
	Since $M = 10^5$, we obtain $\log_K \rho > 0.0002$. Hence, $t(G) \ge n^{2.0002}$ for all sufficiently large $n$.
\end{proof}

\section*{Declaration of generative AI use}

The proof is due to OpenAI's GPT~6.2 Astra. The author's prompt asked the question without suggesting any ideas for a proof. The author subsequently worked through and simplified the generated argument, rewrote it in his own words with feedback from the model, and added further observations, explanations, and references. He takes full responsibility for the contents of this paper.

\end{document}